\documentclass{amsart}
\usepackage{amsmath,amssymb,amsthm,mathtools}
\usepackage{enumitem}
\usepackage{booktabs}
\usepackage{hyperref}
\usepackage{tikz}
\usetikzlibrary{arrows.meta,positioning}
\usetikzlibrary{calc,fit,positioning,arrows.meta}

\newtheorem{theorem}{Theorem}[section]
\newtheorem{lemma}[theorem]{Lemma}
\newtheorem{proposition}[theorem]{Proposition}

\numberwithin{equation}{section}

\theoremstyle{definition}
\newtheorem{definition}[theorem]{Definition}

\theoremstyle{remark}
\newtheorem{remark}[theorem]{Remark}
\newtheorem{example}[theorem]{Example}

\newcommand{\Q}{\mathbb Q}

\newcommand{\cE}{\mathcal E}
\newcommand{\cP}{\mathcal P}
\newcommand{\cG}{\mathcal G}

\newcommand{\ha}{\widehat a}
\newcommand{\Ps}{\mathcal \cE^{\sharp}}
\newcommand{\Ls}{L^{\sharp}}
\newcommand{\Kb}{\overline K}
\newcommand{\Ht}{\widetilde H}

\newcommand{\zero}{\hat 0}
\newcommand{\one}{\hat 1}

\newcommand{\wideand}{\quad \text{ and } \quad}
\DeclareMathOperator{\ex}{ex}
\DeclareMathOperator{\sa}{sa}

\title{The graph $a$-Tur\'an theorem}

\author[S. Choi]{Suyoung Choi}
\address{Department of mathematics, Ajou University, 206, World cup-ro, Yeongtong-gu, Suwon 16499, Republic of Korea}
\email{schoi@ajou.ac.kr}

\author[Y. Yoon]{Younghan Yoon}
\address{Research Institute of Basic Sciences, Ajou University, 206, World cup-ro, Yeongtong-gu, Suwon 16499, Republic of Korea}
\email{younghan300@ajou.ac.kr}
\thanks{Both authors were supported by the National Research Foundation of Korea (NRF) grant funded by the Korean Government (RS-2025-00521982), with additional support for the second named author under grant RS-2021-NR060141.
}

\date{\today}

\begin{document}

\begin{abstract}
We prove analogues of Tur\'an's theorem and the Erd\H{o}s--Stone theorem for graph $a$-numbers.
The $i$th $a$-number~$a_i(G)$ equals the $i$th rational Betti number of the real toric manifold associated with the graph associahedron of~$G$.
We show that the Tur\'an graph~$T_r(n)$ maximizes every $a$-number among~$K_{r+1}$-free graphs on~$n$ vertices.
We also show that for any fixed graph~$H$ with at least one edge and each fixed~$i\geq1$, the limiting extremal density of the $i$th $a$-number among $H$-free graphs depends only on the chromatic number of~$H$.
Since~$a_1(G)$ is the number of edges of~$G$, these results recover the two classical theorems in degree one and extend them to higher $a$-numbers.
\end{abstract}
\maketitle
  \tableofcontents
  
  \section{Introduction}
All graphs in this paper are finite and simple.
For a graph~$G$, we write~$V(G)$ for its vertex set, $e(G)$ for its
number of edges, and~$G[S]$ for the subgraph induced by
$S\subseteq V(G)$.
The graph with no vertices is denoted by~$\varnothing$.
We write~$K_n$ for the complete graph on~$n$ vertices and~$\Kb_n$
for its complement, the graph on~$n$ vertices with no edges.
A graph is \emph{$K_{r+1}$-free} if it has no clique of size~$r+1$.

Let~$n,r\geq1$.
The \emph{Tur\'an graph}~$T_r(n)$ is the complete $r$-partite graph
on~$n$ vertices whose~$r$ part sizes differ by at most one.
Empty parts are allowed, and~$T_r(n)=K_n$ when~$r\geq n$.
The case~$r=2$ is Mantel's theorem~\cite{Mantel1907}.
Tur\'an proved the general result~\cite{Turan1941}.
Every $K_{r+1}$-free graph~$G$ on~$n$ vertices satisfies
$$
e(G)\leq e(T_r(n)),
$$
with equality if and only if~$G\cong T_r(n)$.
Zykov later strengthened this result by showing that~$T_r(n)$
maximizes the number of cliques of each fixed order among
$K_{r+1}$-free graphs on~$n$ vertices~\cite{Zykov1949}.
The Erd\H{o}s--Stone theorem~\cite{Erdos-Stone1946} showed that
Tur\'an graphs also govern the asymptotic edge-extremal problem for
arbitrary forbidden graphs, with the chromatic number determining
the leading term.
For a historical account of the theory developed from
Tur\'an's theorem, see~\cite{Simonovits1977}.

We study analogous extremal problems for graph $a$-numbers.
The \emph{signed $a$-number}~$\sa(G)$ was introduced
in~\cite{Choi-Park2015} and is defined recursively as follows.
\begin{enumerate}[label=(\roman*)]
    \item\label{sa:empty}
    $\sa(\varnothing)=1$.

    \item\label{sa:odd}
    If~$G$ has a connected component of odd order, then~$\sa(G)=0$.

    \item\label{sa:rec}
    If~$G$ is nonempty and every connected component of~$G$ has even
    order, then
    $$
    \sa(G)
    =
    -\sum_{S\subsetneq V(G)}\sa(G[S]).
    $$
\end{enumerate}
The \emph{$a$-number} of~$G$ is~$a(G)=|\sa(G)|$, and for~$i\geq0$
the \emph{$i$th $a$-number} of~$G$ is
$$
a_i(G)
=
\sum_{\substack{S\subseteq V(G)\\ |S|=2i}}
a(G[S]).
$$
In particular,
$$
a_0(G)=1,\quad
a_1(G)=e(G),\wideand
a_i(G)=0\text{ if }2i>|V(G)|.
$$

Let~$M(G)$ be the real toric manifold associated with the graph
associahedron of~$G$.
See~\cite{Seo-Shin2015,Dotsenko-Shadrin-Vallette2019,Park-Park-Park2020,Park-Park2023,Choi-Yoon2023,Choi-Yoon2026,Yoon2025,Choi-Yoon2025} for topological and combinatorial studies of these manifolds and related structures.
Let~$\beta_i(M(G))$ denote its $i$th rational Betti number.
It was proved in~\cite{Choi-Park2015} that
\begin{equation}\label{eq:betti}
\beta_i(M(G))=a_i(G)
\end{equation}
for every graph~$G$.

Our first result is a Tur\'an theorem for graph $a$-numbers.

\begin{theorem}\label{thm:main}
Let~$n,r\geq1$, and let~$G$ be a $K_{r+1}$-free graph on~$n$
vertices.
Then
$$
a_i(G)\leq a_i(T_r(n))
\qquad\text{for every }i\geq0.
$$
If~$1\leq i\leq n/2$, then equality holds if and only if
$G\cong T_r(n)$.
\end{theorem}

For~$i=1$, Theorem~\ref{thm:main} is Tur\'an's theorem.
When~$r\geq n$, we have~$T_r(n)=K_n$, and the inequality also follows
from~\cite{Choi-Yoon2025}, where it was proved that
$a_i(G)\geq a_i(H)$ for every~$i$ whenever~$H$ is a subgraph of~$G$.

In fact, for~$1\leq i\leq n/2$, we prove the stronger inequality
$$
a_i(T_r(n))-a_i(G)
\geq
\frac{e(T_r(n))-e(G)}{n-1},
$$
as stated in Remark~\ref{rem:deficit}.

By~\eqref{eq:betti}, Theorem~\ref{thm:main} gives
$$
\beta_i(M(G))
\leq
\beta_i(M(T_r(n)))
$$
for every $K_{r+1}$-free graph~$G$ on~$n$ vertices and every~$i\geq0$.
Thus, among the real toric manifolds associated with these graphs,
$M(T_r(n))$ maximizes all rational Betti numbers simultaneously.
Equality in any one degree~$1\leq i\leq n/2$ already forces
$G\cong T_r(n)$.

Our second result is an Erd\H{o}s--Stone theorem for graph
$a$-numbers.
Let~$H$ be a fixed graph with at least one edge.
A graph is \emph{$H$-free} if it has no subgraph, not necessarily
induced, isomorphic to~$H$.
For~$i\geq1$, put
$$
\ex_i(n,F)
\coloneq
\max\{a_i(G)\mid |V(G)|=n,\ G\text{ is }H\text{-free}\}.
$$

\begin{theorem}\label{thm:erdos-stone}
Let~$H$ be a fixed graph with at least one edge, let
$r=\chi(H)-1$, and fix~$i\geq1$.
Then, as~$n\to\infty$,
$$
\ex_i(n,F)
=
a_i(T_r(n))+o(n^{2i}).
$$
\end{theorem}

For~$i=1$, Theorem~\ref{thm:erdos-stone} recovers the classical
Erd\H{o}s--Stone asymptotic formula~\cite{Erdos-Stone1946}.
We also give explicit formulas for the $a$-numbers of Tur\'an graphs
in Lemma~\ref{cor:turan-values}.

The paper is organized as follows.
Section~\ref{sec:aux} introduces a homological upper bound for graph
$a$-numbers and establishes its basic properties.
Section~\ref{sec:Zykov} develops the Zykov symmetrization needed for
the extremal problem.
Section~\ref{sec:proof} compares complete multipartite graphs and proves
Theorem~\ref{thm:main}.
Section~\ref{sec:examples} derives explicit formulas for the $a$-numbers
of complete multipartite and Tur\'an graphs, and proves
Theorem~\ref{thm:erdos-stone} together with the density formula.

\section{A homological upper bound for graph $a$-numbers}\label{sec:aux}

Let~$\cP$ be a finite poset, and let $\Delta(\cP)$ be its order complex.
If $\cP$ has a least element $\zero$ and a greatest element $\one$, write $\overline \cP=\cP\setminus\{\zero,\one\}$ for its proper part and~$\mu_{\cP}$ for its M\"obius function.
All reduced homology groups have coefficients in~$\Q$, with $\Ht_{-1}(\{\varnothing\})=\Q$.

A graded poset~$\cP$ is called \emph{locally semimodular} if, in every interval of~$\cP$, any two distinct elements covering a common element have a common cover.
If~$\cP$ is a bounded locally semimodular graded poset of rank~$m\ge1$, then
$\Delta(\overline{\cP})$ is pure of dimension~$m-2$ and has reduced homology only in that degree~\cite{Bjorner1980}. Moreover, it is known~\cite{Stanley1997book} that
\begin{equation}\label{eq:poset-homology}
    \mu_\cP(\zero,\one)=(-1)^m\dim\Ht_{m-2}(\Delta(\overline \cP)).
\end{equation}

Throughout this paper, all graphs are finite and simple.
For a graph~$G$, let~$V(G)$ denote its vertex set and~$e(G)$ its
number of edges. For~$S\subseteq V(G)$, let~$G[S]$ denote the
subgraph induced by~$S$. The graph with no vertices is denoted
by~$\varnothing$.
Let~$\cE(G)$ be the poset of subsets~$S\subseteq V(G)$ such that every component of~$G[S]$ has even
order, ordered by inclusion.

Let~$H$ be a graph with~$2m\ge2$ vertices such that every component
of~$H$ has even order. Then~$\cE(H)$ is graded of rank~$m$, with
least element~$\varnothing$ and greatest element~$V(H)$.
In addition,
\begin{equation}\label{eq:a-homology}
\dim\Ht_k\bigl(\Delta(\overline{\cE(H)})\bigr)
=
\begin{cases}
a(H), & k=m-2,\\
0, & k\ne m-2.
\end{cases}
\end{equation}
Equation~\eqref{eq:a-homology} was proved for connected graphs
in~\cite{Choi-Park2015}, but it remains valid without the
connectedness assumption.

We now replace the even-component condition by the weaker requirement that an even set contain an edge.
Assume that~$H$ is a graph with $2m\ge2$ vertices and at least one edge.
We define
$$
    \Ps(H)=\{\varnothing\}\cup\bigl\{S\subseteq V(H) \mid |S|\text{ is even and $H[S]$ has an edge}\bigr\},
$$
ordered by inclusion.
Note that the vertices of the order complex~$\Delta\bigl(\overline{\Ps(H)}\bigr)$ are the nonempty proper subsets $S\subsetneq V(H)$ of even size such that $H[S]$ has an edge, and the faces are the chains of such sets.

\begin{proposition}\label{lem:sharp}
    A poset~$\Ps(H)$ is graded of rank $m$ and locally semimodular.
    Consequently, $\Delta\bigl(\overline{\Ps(H)}\bigr)$ has reduced homology only in degree $m-2$.
\end{proposition}
\begin{proof}
    The least element of $\Ps(H)$ is $\varnothing$, and since $H$ has an edge, $V(H)$ is the greatest element.
    Let $S\in\Ps(H)$ be nonempty.
    If $S\subseteq T\subseteq V(H)$ and $|T|$ is even, then $T\in\Ps(H)$, because $T$ contains an edge of $H[S]$.
    So the elements covering $S$ are the sets $S\cup\{x,y\}$ with $x,y\notin S$.
    The elements covering $\varnothing$ are the edges of $H$, and every nonempty element of $\Ps(H)$ contains one of them.
    Hence $\Ps(H)$ is graded, with rank function $S\mapsto|S|/2$.
    
    Now let $U\ne U'$ both cover $T$, and let $U,U'\subseteq B$ with $B\in\Ps(H)$.
    Then $|U|=|U'|=|T|+2$, so $|U\cup U'|$ is $|T|+3$ or $|T|+4$.
    If $|U\cup U'|=|T|+4$, let $Z=U\cup U'$.
    If $|U\cup U'|=|T|+3$, then $|U\cup U'|$ is odd while $|B|$ is even, so there is $q\in B\setminus(U\cup U')$; let $Z=U\cup U'\cup\{q\}$.
    In both cases $Z\subseteq B$, $|Z|=|U|+2$, and $Z$ contains an edge, because $U$ does.
    So $Z\in\Ps(H)$, and $Z$ covers $U$ and $U'$.
    The homology statement follows from~\eqref{eq:poset-homology}.
\end{proof}

\begin{definition}\label{def:ha}
    For a graph $H$ on $2m$ vertices, let
    $$
        \ha(H)= \begin{cases}
            1 & \text{if } H=\varnothing,\\
            0 & \text{if } H\ne\varnothing \text{ has no edges},\\
            \dim\Ht_{m-2}(\Delta\bigl(\overline{\Ps(H)}\bigr)) & \text{if } m \geq 1 \text{ and $H$ has an edge}.
        \end{cases}
    $$
    For a graph $G$ and $i\ge0$, let
    $$
        \ha_i(G)=\sum_{\substack{S\subseteq V(G)\\ |S|=2i}}\ha(G[S]).
    $$
\end{definition}

Clearly, $\ha_i(G)\ge0$ for each~$i\ge0$.
In particular, $\ha_0(G)=1$
and~$\ha_i(G)=0$ for~$2i>|V(G)|$.
Since $\ha(K_2)=1$ and pairs of nonadjacent vertices contribute zero, we also have $\ha_1(G)=e(G)$.

We write~$P_n$ for the path graph on~$n$ vertices.
\begin{example}\label{ex:ha-four}
Let~$H$ be a graph on four vertices with at least one edge.
The poset~$\overline{\Ps(H)}$ is a rank-zero poset whose elements are precisely the edges of~$H$.
Hence, its order complex consists of~$e(H)$ points, and therefore
$$
\ha(H)=e(H)-1.
$$
For instance, $\ha(P_3\sqcup K_1)=1$, while
$a(P_3\sqcup K_1)=0$.
\end{example}

A graph is \emph{complete multipartite} if its vertex set can be partitioned into parts such that two vertices are adjacent exactly
when they lie in different parts.

\begin{proposition}\label{prop:compare}
    For every~$i\ge0$, we have
$$
a_i(G)\le \ha_i(G),
$$
with equality for all~$i$ whenever~$G$ is complete multipartite.
\end{proposition}
\begin{proof}
    Induced subgraphs of complete multipartite graphs are complete multipartite.
    Then, it is enough to show that $a(H)\le\ha(H)$ for every graph $H$ with an even number of vertices, with equality when $H$ is complete multipartite.
    
    For convenience, write
$$
L(H)=\Delta\bigl(\overline{\cE(H)}\bigr)
\qquad\text{and}\qquad
\Ls(H)=\Delta\bigl(\overline{\Ps(H)}\bigr).
$$
    
    If~$H=\varnothing$, then~$a(H)=\ha(H)=1$.
    If~$H$ has a component of odd order, then~$a(H)=0\le\ha(H)$.
    Suppose now that~$|V(H)|=2m\ge2$ and every component of $H$ has even order.
    Then~$H$ has an edge.
    Moreover, every nonempty~$S\in\cE(H)$ contains an edge, because each component of~$H[S]$ is connected and has at least two vertices.
    It follows that~$\cE(H)\subseteq\Ps(H)$, and hence~$L(H)$ is a subcomplex of~$\Ls(H)$.
    Put~$d=m-2$.
    Neither complex has a face of dimension $d+1$, so for each of them,~$\Ht_d$ is the space of $d$-cycles.
    A $d$-cycle of~$L(H)$ is also a $d$-cycle of~$\Ls(H)$.
    Hence, by~\eqref{eq:a-homology},
    $$
        a(H)=\dim\Ht_d(L(H))\le\dim\Ht_d(\Ls(H))=\ha(H).
    $$
    
    Now let $H$ be complete multipartite with $|V(H)|$ even.
    If $H$ is nonempty and has no edges, then all components of $H$ are single vertices, so $a(H)=0=\ha(H)$.
    Suppose that $H$ has an edge.
    If $S\subseteq V(H)$ meets at least two parts, then $H[S]$ is connected: two vertices of $S$ in different parts are adjacent, and two vertices of $S$ in the same part have a common neighbor in $S$.
    In particular, $H$ is connected.
    For a nonempty set $S\subseteq V(H)$ of even size, this gives
    $$
        H[S]\text{ has an edge}  \iff S\text{ meets at least two parts}  \iff S\in\cE(H).
    $$
    Thus, $\cE(H)=\Ps(H)$ and~$L(H)=\Ls(H)$.
    Therefore,~\eqref{eq:a-homology} gives $a(H)=\ha(H)$.
\end{proof}

The invariant~$\ha_i(G)$ admits a formula in terms of independent sets.
Let $I_k(G)$ be the number of independent sets of size $k$ in $G$ with $I_0(G)=1$.
Let $E_{2i}$ denote the Euler secant numbers~\cite[A000364]{oeis}. 

\begin{theorem}\label{prop:independent}
    Let $G$ be a graph on $n$ vertices and let $0\leq i \leq \lfloor \frac{n}{2} \rfloor$.
    Then
    $$
    \ha_i(G)=\sum_{j=0}^i(-1)^j E_{2i-2j} \binom{n-2j}{2i-2j}I_{2j}(G).
    $$
\end{theorem}
\begin{proof}
    First let $H$ have $2i$ vertices.
    For each even set $S\subseteq V(H)$, put
    $$
        b(S)=(-1)^{\frac{|S|}{2}}\ha(H[S]).
    $$
    If $H[S]$ has an edge, Proposition~\ref{lem:sharp} and \eqref{eq:poset-homology} identify $b(S)$ with the M\"obius number of $\Ps(H[S])$.
    Its defining recurrence, together with the definition for edgeless graphs, gives
    $$
        \sum_{\substack{T\subseteq S\\|T|\text{ even}}}b(T) =\begin{cases}
            1 & \text{if $S$ is independent},\\
            0 & \text{otherwise}.
        \end{cases}
    $$
    In the poset of all even subsets of $V(H)$, the M\"obius function is
    $$
        \mu(T,S)=(-1)^{(|S|-|T|)/2}E_{|S|-|T|}.
    $$
    This follows by induction from the coefficient identity $\cos x\sec x=1$.
    M\"obius inversion therefore gives
    \begin{equation}\label{eq:independent-local}
        \ha(H)=\sum_{\substack{T\subseteq V(H)\\T\text{ independent}\\|T|\text{ even}}} (-1)^{|T|/2}E_{2i-|T|}.
    \end{equation}
    This also holds when $H$ is empty.
    Now sum \eqref{eq:independent-local} over all $2i$-subsets of $V(G)$.
    Each independent set of size $2j$ is contained in exactly $\binom{n-2j}{2i-2j}$ such subsets, which proves the formula.
\end{proof}

\section{Zykov symmetrization}\label{sec:Zykov}

Zykov's symmetrization method~\cite{Zykov1949} is a classical tool in extremal graph theory for reducing problems to complete multipartite graphs.
In this section, we show that~$\ha_i$ satisfies a natural inequality under Zykov symmetrization.

For a vertex~$u$ of a graph~$G$, let~$N_G(u)$ denote its neighborhood.
For distinct nonadjacent vertices~$u,v$, let~$G+uv$ denote the graph
obtained from~$G$ by adding the edge~$uv$.
For a finite set~$U$, let~$K_U$ denote the complete graph on~$U$.
For graphs~$G$ and~$H$ on disjoint vertex sets, their
\emph{join}~$G\vee H$ is obtained from~$G\sqcup H$ by adding all
edges between~$V(G)$ and~$V(H)$.

To begin with, we describe how~$\ha_i(G)$ changes when an edge is added.
Let~$u$ and~$v$ be distinct nonadjacent vertices of~$G$, and set
$$
W=V(G)\setminus\{u,v\}
\wideand
C=W\setminus\bigl(N_G(u)\cup N_G(v)\bigr).
$$
Let~$R_{uv}(G)$ be the graph on~$W$ obtained from~$G[W]$ by
joining every vertex of~$W\setminus C$ to every other vertex of~$W$.
Equivalently,
$$
R_{uv}(G)=G[C]\vee K_{W\setminus C}.
$$
For instance, Figure~\ref{fig:Ruv} illustrates the construction.

\begin{figure}
\centering
\begin{tikzpicture}[
    >=Latex,
    vertex/.style={circle, draw, thick, minimum size=18pt, inner sep=1pt, font=\small},
    lab/.style={draw=none, font=\small},
    title/.style={draw=none, font=\normalsize},
    region/.style={draw, dashed, rounded corners, inner sep=6pt}
]

\begin{scope}[xshift=0cm]
\node[title] at (2,-3) {$G$};

\node[vertex] (u) at (4,1.0) {$u$};
\node[vertex] (v) at (4,-1.0) {$v$};

\node[vertex] (a) at (0.3,1.4) {$a$};
\node[vertex] (b) at (0.3,-0.1) {$b$};

\node[vertex] (c) at (2.5,1.6) {$c$};
\node[vertex] (d) at (1.7,0.0) {$d$};
\node[vertex] (e) at (2.5,-1.6) {$e$};

\draw[thick] (a)--(b);
\draw[thick] (b)--(c);
\draw[thick] (c)--(d);
\draw[thick] (c)--(e);

\draw[thick] (u)--(c);
\draw[thick] (u)--(d);
\draw[thick] (u)--(e);

\draw[thick] (v)--(c);
\draw[thick] (v)--(e);

\node[region, fit=(a)(b), label=above:{\small $C$}] {};
\node[region, fit=(c)(d)(e), label=above:{\small $W\setminus C$}] {};
\end{scope}

\node[draw=none, font=\Large] at (5.7,0) {$\Longrightarrow$};

\begin{scope}[xshift=7cm]
\node[title] at (1.5,-3) {$R_{uv}(G)$};

\node[vertex] (a2) at (0.3,1.4) {$a$};
\node[vertex] (b2) at (0.3,-0.1) {$b$};

\node[vertex] (c2) at (2.5,1.6) {$c$};
\node[vertex] (d2) at (1.7,0.3) {$d$};
\node[vertex] (e2) at (2.5,-1.6) {$e$};

\draw[thick] (a2)--(b2);

\draw[thick] (c2)--(d2);
\draw[thick] (c2)--(e2);
\draw[thick] (d2)--(e2);

\draw[thick] (a2)--(c2);
\draw[thick] (a2)--(d2);
\draw[thick] (a2)--(e2);
\draw[thick] (b2)--(c2);
\draw[thick] (b2)--(d2);
\draw[thick] (b2)--(e2);

\node[region, fit=(a2)(b2), label=above:{\small $C$}] {};
\node[region, fit=(c2)(d2)(e2), label=above:{\small $W\setminus C$}] {};
\end{scope}
\end{tikzpicture}
\caption{$R_{uv}(G)$}
\label{fig:Ruv}
\end{figure}

\begin{proposition}\label{lem:positive}
    Let $H$ be a graph with $2m\ge2$ vertices.
    \begin{enumerate}
      \item\label{lem:positive1} For any distinct nonadjacent vertices~$u$ and~$v$ of~$H$,
          $$
          \ha(H+uv)=\ha(H)+\ha(R_{uv}(H)).
          $$
      \item\label{lem:positive2} If~$H$ has a universal vertex, then~$\ha(H)\ge1$.
    \end{enumerate}
\end{proposition}

\begin{proof}
We first prove~\eqref{lem:positive1}, and then use it to prove~\eqref{lem:positive2}.
\begin{enumerate}
  \item Let $W$ and $C$ be defined as above, and let $R:=R_{uv}(H)$.
    Adding $uv$ removes exactly the independent sets of the form $\{u,v\}\cup T$, where $T$ is independent in~$H[C]$.
    Thus Theorem~\ref{prop:independent}, applied to~$H+uv$ and~$H$ with~$i=m$, gives
    $$
        \ha(H+uv)-\ha(H)=\sum_{\substack{T\subseteq C\\T\text{ independent in }H[C]\\|T|\text{ even}}} (-1)^{|T|/2}E_{2m-2-|T|}.
    $$
    Every vertex of $W\setminus C$ is adjacent to all other vertices of $R$.
    Hence the independent sets of even size in $R$ are exactly those in $H[C]$, including the empty set.
    Applying Theorem~\ref{prop:independent} to~$R$ with~$i=m-1$ shows that the sum is~$\ha(R)$.
  \item
  Let $s$ be a universal vertex of~$H$.
  We use induction on~$m$.
    If~$m=1$, then~$H=K_2$ and~$\ha(H)=1$.
    
    Suppose~$m\ge2$.
    Choose a vertex~$x\neq s$, and let~$H_x$ be obtained from~$H$ by deleting the edge~$sx$.
    The vertex $s$ is adjacent to every vertex of~$W:=V(H)\setminus\{s,x\}$.
    Hence~$s$ and~$x$ have no common non-neighbor in~$W$, which implies that~$R_{sx}(H_x)=K_W$.
    Combining~\eqref{lem:positive1} with the induction hypothesis applied to~$K_W$, we obtain
    $$
        \ha(H)=\ha(H_x)+\ha(K_W)\ge\ha(K_W)\ge1,
    $$
    as desired.
\end{enumerate}
\end{proof}

If $I\subseteq V(G)$ contains $u$ and $v$, then the common non-neighbors of $u$ and $v$ in~$G[I]$ are the vertices of~$C\cap I$.
Hence
\begin{equation}\label{eq:R-restrict}
    R_{uv}(G[I])=R_{uv}(G)[I\setminus\{u,v\}].
\end{equation}
\begin{lemma}\label{prop:edge}
    Let $u,v$ be nonadjacent vertices of a graph $G$.
    Then
    $$
        \ha_i(G+uv)-\ha_i(G)=\ha_{i-1}(R_{uv}(G)) \quad \text{ for every }i\ge1.
    $$
\end{lemma}
\begin{proof}
    Let $I\subseteq V(G)$ with $|I|=2i$.
    If $I$ does not contain both $u$ and $v$, then $(G+uv)[I]=G[I]$.
    If $u,v\in I$, then $(G+uv)[I]=G[I]+uv$, and Proposition~\ref{lem:positive}\eqref{lem:positive1} and~\eqref{eq:R-restrict} give
    $$
        \ha\bigl((G+uv)[I]\bigr)-\ha(G[I])=\ha\bigl(R_{uv}(G)[I\setminus\{u,v\}]\bigr).
    $$
    As $I$ runs over the $2i$-subsets of $V(G)$ that contain $u$ and $v$, the set $I\setminus\{u,v\}$ runs over the $(2i-2)$-subsets of $V(R_{uv}(G))$.
    Summing proves the identity.
    The last statement follows because $\ha_{i-1}\ge0$ and $\ha_0=1$.
\end{proof}
Consequently,~$\ha_i$ is monotone with respect to edge inclusion.

\begin{definition}\label{def:copy}
    Let $u,v$ be nonadjacent vertices of $G$.
    The graph $G_{u\leftarrow v}$ is obtained from $G$ by deleting all edges at $u$ and then joining $u$ to every vertex of $N_G(v)$.
\end{definition}

The vertices~$u$ and~$v$ are nonadjacent in~$G_{u\leftarrow v}$ and both have neighborhood~$N_G(v)$.
All edges not incident to~$u$ are unchanged.
Moreover,
$$
    e(G_{u\leftarrow v})+e(G_{v\leftarrow u})=2e(G).
$$

\begin{proposition}\label{lem:clique}
    If $G$ is $K_{r+1}$-free, then so is $G_{u\leftarrow v}$.
\end{proposition}
\begin{proof}
    Let $Q$ be a clique of $G_{u\leftarrow v}$.
    If $u\notin Q$, then $Q$ is a clique of $G$.
    If $u\in Q$, then $v\notin Q$ and $Q\setminus\{u\}\subseteq N_G(v)$.
    Therefore, $(Q\setminus\{u\})\cup\{v\}$ is a clique of $G$ of the same size.
\end{proof}

\begin{lemma}\label{lem:sym-edge}
    If~$u$ and~$v$ are nonadjacent vertices of $G$, then
    $$
        \ha_i(G_{u\leftarrow v}+uv)+\ha_i(G_{v\leftarrow u}+uv)=2\ha_i(G+uv) \qquad \text{for each } i \geq 0.
    $$
\end{lemma}
\begin{proof}
    If~$2i>|V(G)|$, then the desired identity is immediate, since all three terms are zero.
Thus, we may assume that~$2i\le |V(G)|$.
By Theorem~\ref{prop:independent}, it is enough to show that
$$
I_k(G_{u\leftarrow v}+uv)+I_k(G_{v\leftarrow u}+uv)
=
2I_k(G+uv)
$$
for every~$k\geq0$.

Since $uv$ is an edge, no independent set in these graphs contains both~$u$ and~$v$.
    Independent sets containing neither vertex are the same in all three graphs.
    For the sets containing exactly one of them, put $W=V(G)\setminus\{u,v\}$.
    Let $p$ and $q$ be the numbers of independent $(k-1)$-sets in $G[W\setminus N_G(u)]$ and $G[W\setminus N_G(v)]$, respectively.
    The contributions of these independent sets to
$$
I_k(G_{u\leftarrow v}+uv),\qquad
I_k(G_{v\leftarrow u}+uv),\qquad
I_k(G+uv)
$$
are $2q$, $2p$, and $p+q$, respectively.
    This proves the identity for $k\ge1$, while $k=0$ is immediate.
\end{proof}

Zykov symmetrization is particularly useful when one of the two
symmetrizations can be chosen without decreasing the graph parameter.
We have the following stronger inequality for ~$\ha_i$.

\begin{theorem}\label{prop:sym}
    Let $u,v$ be nonadjacent vertices of $G$.
    For every $i\geq0$,
    $$
        \ha_i(G_{u\leftarrow v})+\ha_i(G_{v\leftarrow u})\ge2\ha_i(G).
    $$
\end{theorem}

\begin{proof}
For~$i=0$, both sides equal~$2$.
Suppose~$i\ge1$, and write~$G_1=G_{u\leftarrow v}$ and
$G_2=G_{v\leftarrow u}$.
By Lemma~\ref{prop:edge}, for~$\cG \in \{G ,G_1, G_2\}$,
$$
\ha_i(\cG)
=
\ha_i(\cG+uv)-\ha_{i-1}(R_{uv}(\cG)).
$$
Hence, using Lemma~\ref{lem:sym-edge},
\begin{align*}
  \ha_i(G_1)+\ha_i(G_2)-2\ha_i(G)
&=
\ha_i(G_1+uv)+\ha_i(G_2+uv)-2\ha_i(G+uv)
\notag\\
&\quad
-\ha_{i-1}(R_{uv}(G_1))
-\ha_{i-1}(R_{uv}(G_2))
+2\ha_{i-1}(R_{uv}(G))
\notag\\
&=
2\ha_{i-1}(R_{uv}(G))
-\ha_{i-1}(R_{uv}(G_1))
-\ha_{i-1}(R_{uv}(G_2)).
  \end{align*}

To finish the proof, it is sufficient to show that
\begin{equation}\label{eq:pro}
  \ha_{i-1}(R_{uv}(G))
-\ha_{i-1}(R_{uv}(G_j)) \geq 0 \qquad \text{for } j=1,2.
\end{equation}
Let~$W=V(G)\setminus\{u,v\}$.
Since only edges incident to~$u$ or~$v$ are changed,
$$
G[W]=G_1[W]=G_2[W].
$$
Therefore,
\begin{align*}
R_{uv}(G_1)
&=G[W\setminus N_G(v)]\vee K_{N_G(v)}\\
&\subseteq G[W\setminus(N_G(u)\cup N_G(v))]
\vee K_{N_G(u)\cup N_G(v)}=R_{uv}(G),\\
R_{uv}(G_2)
&=G[W\setminus N_G(u)]\vee K_{N_G(u)}\\
&\subseteq G[W\setminus(N_G(u)\cup N_G(v))]
\vee K_{N_G(u)\cup N_G(v)}=R_{uv}(G).
\end{align*}
Since~$\ha_{i-1}$ is monotone with respect to edge inclusion,
$$
\ha_{i-1}(R_{uv}(G_j))
\le\ha_{i-1}(R_{uv}(G))
\qquad\text{for }j=1,2.
$$
Thus~\eqref{eq:pro} holds, completing the proof.
\end{proof}

We write~$C_n$ for the cycle graph on~$n$ vertices.
\begin{remark}
For the graph~$a$-numbers, the inequality in Theorem~\ref{prop:sym} does not hold in general.
For example, let~$G=P_4$ be the path on four vertices labeled~$1,2,3,4$ in this order.
Then
$$
G_{1\leftarrow3}\cong C_4
\qquad\text{and}\qquad
G_{3\leftarrow1}\cong P_3\sqcup K_1.
$$
Since
$$
a_2(P_4)=2,\qquad
a_2(C_4)=3,\qquad
a_2(P_3\sqcup K_1)=0,
$$
we have
$$
a_2(G_{1\leftarrow3})+a_2(G_{3\leftarrow1})
=3<4=2a_2(G).
$$
\end{remark}

\section{Proof of Theorem~\ref{thm:main}}\label{sec:proof}

To prove Theorem~\ref{thm:main}, we compare complete multipartite graphs by moving one vertex from a larger part to a smaller one.
For integers~$0\le z\le N$, let
$$
    S_{N,z}=\Kb_z\vee K_{N-z},
$$
where~$\Kb_z$ is the discrete graph on~$z$ vertices and~$K_{N-z}$ is the complete graph on~$N-z$ vertices.
This is the complete multipartite graph on~$N$ vertices with one part of size~$z$ and~$N-z$ parts of size one.
For disjoint sets~$V_1,\ldots,V_r$, we write
$K(V_1,\ldots,V_r)$ for the complete multipartite graph with parts
$V_1,\ldots,V_r$, allowing empty parts.

\begin{lemma}\label{lem:mo}
For the complete multipartite graph~$K=K(V_1,\ldots,V_r)$ with~$r\ge2$, the following hold.
\begin{enumerate}
    \item\label{lem:mo1}
    $e(\widetilde{K})-e(K)=|V_1|-|V_2|-1.$

    \item\label{lem:mo2} For each~$i\geq1$,
    $$
\ha_i(\widetilde{K})-\ha_i(K)
=
\sum_{\ell=|V_2|}^{|V_1|-2}
\ha_{i-1}(S_{n-2,\ell}).
$$

    \item\label{lem:mo3} 
Assume that~$|V_1| \geq |V_2|+2$.
For~$1\le i\le n/2$,
$$
\ha_i(\widetilde{K})-\ha_i(K)\ge1,
$$
where~$\widetilde{K}
=
K(V_1\setminus\{w\},V_2\cup\{w\},V_3,\ldots,V_r)$ for a fixed~$w\in V_1$.
\end{enumerate}
\end{lemma}
\begin{proof}
    The edges between $V_1\cup V_2$ and the other parts do not change.
    The number of edges within $V_1\cup V_2$ changes from
$|V_1||V_2|$ to $(|V_1|-1)(|V_2|+1)$.
Hence
$$
e(\widetilde{K})-e(K)
=
(|V_1|-1)(|V_2|+1)-|V_1||V_2|
=
|V_1|-|V_2|-1,
$$
which proves~\eqref{lem:mo1}.

Let~$G_0$ be the complete multipartite graph with parts
$\{w\}$, $V_1\setminus\{w\}$, $V_2,\ldots,V_r$.
Then~$K$ is obtained from~$G_0$ by deleting all edges between~$w$ and~$V_1\setminus\{w\}$, and~$\widetilde{K}$ is obtained from~$G_0$ by
deleting all edges between~$w$ and $V_2$.

Let~$U=\{u_1,\ldots,u_k\}$ be a part of~$G_0$ other than $\{w\}$.
For~$1\leq t\leq k$, let~$G_t$ be obtained from~$G_0$
by deleting the edges $wu_1,\ldots,wu_t$.
Then $G_{t-1}=G_t+wu_t$.
The common non-neighbors of~$w$ and~$u_t$ in~$G_t$ are
$u_1,\ldots,u_{t-1}$, which are pairwise nonadjacent.
Hence
$$
R_{wu_t}(G_t)
\cong
\Kb_{t-1}\vee K_{n-t-1}
=
S_{n-2,t-1}.
$$
Applying Lemma~\ref{prop:edge} at each step and summing gives
$$
\ha_i(G)-\ha_i(G_k)
=
\sum_{\ell=0}^{k-1}\ha_{i-1}(S_{n-2,\ell}).
$$
Applying this identity to $U=V_1\setminus\{w\}$ and $U=V_2$, and
subtracting the resulting identities, gives
$$
\ha_i(\widetilde{K})-\ha_i(K)
=
\sum_{\ell=|V_2|}^{|V_1|-2}
\ha_{i-1}(S_{n-2,\ell}),
$$
which proves~\eqref{lem:mo2}.

Suppose that $1\leq i\leq n/2$.
Since all terms in the sum in~\eqref{lem:mo2} are nonnegative,
it suffices to show that
$\ha_{i-1}(S_{n-2,|V_2|})\geq1$.
If $i=1$, then $\ha_0(S_{n-2,|V_2|})=1$.
Assume that $i\geq2$, and then $n\geq4$.
Since $|V_1|\geq|V_2|+2$ and $|V_1|+|V_2|\leq n$, we have
$$
|V_2|\leq\frac{n}{2}-1\leq n-3.
$$
Thus $S_{n-2,|V_2|}$ has a singleton part $\{s\}$, and $s$ is
adjacent to every other vertex.
Since $n\geq2i$, there is a set $T$ of $2i-2$ vertices of
$S_{n-2,|V_2|}$ containing $s$.
By the definition of $\ha_{i-1}$ and
Proposition~\ref{lem:positive}\eqref{lem:positive2},
$$
\ha_{i-1}(S_{n-2,|V_2|})
\geq
\ha\bigl(S_{n-2,|V_2|}[T]\bigr)
\geq1,
$$
which proves~\eqref{lem:mo3}.
\end{proof}

For $i\geq0$ and $\lambda>0$, define
$$
\Phi_{i,\lambda}(G)
=
\ha_i(G)-\frac{e(G)}{\lambda}.
$$

\begin{lemma}\label{prop:multipartite}
Let $n,r\geq1$, $i\geq0$, and $\lambda>0$, and fix a set~$V$ of
$n$ vertices.
Among all $K_{r+1}$-free graphs on~$V$, the function
$\Phi_{i,\lambda}$ attains its maximum at a complete multipartite
graph with at most $r$ nonempty parts.
\end{lemma}

\begin{proof}
For a graph~$H$ on~$V$ and a vertex~$w\in V$, let
$d_H(w)=|N_H(w)|$, and put
$$
D(H)=\sum_{w\in V}d_H(w)^2.
$$
There are finitely many graphs on~$V$.
Among the $K_{r+1}$-free graphs on~$V$ that maximize~$\Phi_{i,\lambda}$,
choose one, say~$G$, for which
$$
D(G)=\sum_{w\in V}d_G(w)^2
$$
is largest.

Let~$u,v$ be distinct nonadjacent vertices of~$G$.
Since
$$
e(G_{u\leftarrow v})+e(G_{v\leftarrow u})=2e(G),
$$
Theorem~\ref{prop:sym} gives
\begin{equation}\label{eq:Phi-convex}
\Phi_{i,\lambda}(G_{u\leftarrow v})
+
\Phi_{i,\lambda}(G_{v\leftarrow u})
\geq
2\Phi_{i,\lambda}(G).
\end{equation}
By Proposition~\ref{lem:clique}, the graphs~$G_{u\leftarrow v}$ and
$G_{v\leftarrow u}$ are $K_{r+1}$-free.
The maximality of~$\Phi_{i,\lambda}(G)$ therefore gives
$$
\Phi_{i,\lambda}(G_{u\leftarrow v})
\leq
\Phi_{i,\lambda}(G)
\quad\text{and}\quad
\Phi_{i,\lambda}(G_{v\leftarrow u})
\leq
\Phi_{i,\lambda}(G).
$$
Together with~\eqref{eq:Phi-convex}, this shows that both
$G_{u\leftarrow v}$ and~$G_{v\leftarrow u}$ maximize~$\Phi_{i,\lambda}$.

We compute
$D(G_{u\leftarrow v})+D(G_{v\leftarrow u})-2D(G)$.
In~$G_{u\leftarrow v}$, both~$u$ and~$v$ have degree~$d_G(v)$,
while in~$G_{v\leftarrow u}$, both have degree~$d_G(u)$.
Thus~$u$ and~$v$ contribute
$$
2d_G(v)^2+2d_G(u)^2
-2\bigl(d_G(u)^2+d_G(v)^2\bigr)
=0.
$$
Let~$w\in V\setminus\{u,v\}$ and~$d=d_G(w)$.
If~$w\in N_G(u)\triangle N_G(v)$, then~$w$ has degree~$d+1$ in one
of~$G_{u\leftarrow v},G_{v\leftarrow u}$ and degree~$d-1$ in the other.
Thus, it contributes
$$
(d+1)^2+(d-1)^2-2d^2=2.
$$
Otherwise,~$w$ has the same degree in~$G$, $G_{u\leftarrow v}$, and
$G_{v\leftarrow u}$.
Hence
$$
D(G_{u\leftarrow v})+D(G_{v\leftarrow u})-2D(G)
=
2|N_G(u)\triangle N_G(v)|.
$$
Since~$G_{u\leftarrow v}$ and~$G_{v\leftarrow u}$ also maximize
$\Phi_{i,\lambda}$, the choice of~$G$ gives
$$
D(G_{u\leftarrow v})\leq D(G)
\quad\text{and}\quad
D(G_{v\leftarrow u})\leq D(G).
$$
Therefore~$N_G(u)=N_G(v)$.

Thus any two distinct nonadjacent vertices of~$G$ have the same
neighborhood.
It follows that nonadjacency is transitive.
Indeed, if~$u$ and~$v$ are nonadjacent, $v$ and~$w$ are nonadjacent,
and~$u\neq w$, then
$$
w\notin N_G(v)=N_G(u).
$$
Thus~$u$ and~$w$ are nonadjacent.
Hence the relation in which two vertices are equal or nonadjacent is an equivalence relation.
Its equivalence classes are independent sets, and any two vertices
in different classes are adjacent.
Thus~$G$ is complete multipartite.
Choosing one vertex from each part gives a clique.
Since~$G$ is $K_{r+1}$-free, it follows that~$G$ has at most~$r$
nonempty parts.
\end{proof}

We now prove the quantitative bound and complete the proof of Theorem~\ref{thm:main}.
\begin{lemma}\label{prop:linear}
    Let $n,r\ge1$ and $1\le i\le n/2$, and let $G$ be a $K_{r+1}$-free graph on $n$ vertices.
    Then
    $$
        \ha_i(T_r(n))-\ha_i(G)\ge\frac{e(T_r(n))-e(G)}{n-1}.
    $$
\end{lemma}

\begin{proof}
Since~$1\leq i\leq n/2$, we have~$n\geq2$.
Let~$V=V(G)$.
By Lemma~\ref{prop:multipartite}, there is a complete multipartite
graph~$K$ on~$V$ with at most~$r$ nonempty parts such that
$$
\Phi_{i,n-1}(G)\leq\Phi_{i,n-1}(K).
$$
Write~$K=K(V_1,\ldots,V_r)$, allowing empty parts.

Suppose that two parts, say~$V_p$ and~$V_q$, satisfy
$|V_p|\geq|V_q|+2$.
Move one vertex from~$V_p$ to~$V_q$.
By Lemma~\ref{lem:mo}\eqref{lem:mo1} and~\eqref{lem:mo3}, this move
increases~$\ha_i$ by at least~$1$ and increases the number of edges by
$$
|V_p|-|V_q|-1\leq n-1.
$$
Hence~$\Phi_{i,n-1}$ does not decrease.
Moreover, the move decreases
$$
\sum_{j=1}^r|V_j|^2
$$
by~$2(|V_p|-|V_q|-1)>0$.
Since this sum is a nonnegative integer, the process terminates after
finitely many moves.

At termination, all part sizes differ by at most one, and the resulting
graph is a copy of~$T_r(n)$.
Therefore
$$
\Phi_{i,n-1}(G)
\leq
\Phi_{i,n-1}(K)
\leq
\Phi_{i,n-1}(T_r(n)).
$$
Equivalently,
$$
\ha_i(T_r(n))-\ha_i(G)
\geq
\frac{e(T_r(n))-e(G)}{n-1}.
$$
\end{proof}

We now prove Theorem~\ref{thm:main}.
\begin{proof}[Proof of Theorem~\ref{thm:main}]
    If $i=0$, both sides equal $1$.
    If $2i>n$, both sides equal $0$.
    Let $1\le i\le n/2$.
    By Tur\'an's theorem, $e(G)\le e(T_r(n))$.
    Since $T_r(n)$ is complete multipartite, Proposition~\ref{prop:compare} and Lemma~\ref{prop:linear} give
    $$
        a_i(G)\le\ha_i(G)\le\ha_i(T_r(n))-\frac{e(T_r(n))-e(G)}{n-1}\le\ha_i(T_r(n))=a_i(T_r(n)).
    $$
    If $a_i(G)=a_i(T_r(n))$, then all these inequalities are equalities.
    In particular, $e(G)=e(T_r(n))$, so $G\cong T_r(n)$ by the equality case of Tur\'an's theorem.
    The converse is clear.
\end{proof}

\begin{remark}\label{rem:deficit}
Let $n,r\geq1$ and $1\leq i\leq n/2$, and let~$G$ be a
$K_{r+1}$-free graph on~$n$ vertices.
Proposition~\ref{prop:compare} and Lemma~\ref{prop:linear} give
$$
a_i(T_r(n))-a_i(G)
\geq
\frac{e(T_r(n))-e(G)}{n-1}.
$$
Indeed, $a_i(G)\leq\ha_i(G)$ and
$a_i(T_r(n))=\ha_i(T_r(n))$.
\end{remark}

\section{Proof of Theorem~\ref{thm:erdos-stone}}\label{sec:examples}

When only the part sizes matter, we write~$K_{n_1,\ldots,n_r}$
for the complete multipartite graph with part sizes
$n_1,\ldots,n_r$.
\begin{example}\label{ex:six}
The graph $T_2(6)$ is $K_{3,3}$.
Its first three $a$-numbers are
$$
a_1(K_{3,3})=9,\qquad a_2(K_{3,3})=39,\qquad a_3(K_{3,3})=31.
$$
To compute $a_2$, note that a four vertex subset meets the two parts in sizes $(1,3)$, $(2,2)$ or $(3,1)$.
Since $a(K_{1,3})=2$ and $a(K_{2,2})=3$, we obtain
$$
a_2(K_{3,3})=3\cdot2+9\cdot3+3\cdot2=39.
$$
The defining recurrence then gives
$$
\sa(K_{3,3})=-(1-9+39)=-31.
$$
Thus every triangle-free graph $G$ on six vertices satisfies
$$
a_1(G)\le9,\qquad a_2(G)\le39,\qquad a_3(G)\le31,
$$
and equality in any one of these inequalities holds if and only if $G\cong K_{3,3}$.
Equivalently, the rational Poincar\'e polynomial of $M(T_2(6))$ is
$$
\sum_{i\ge0}\beta_i(M(T_2(6)))t^i=1+9t+39t^2+31t^3,
$$
and it bounds the rational Poincar\'e polynomials of $M(G)$ coefficientwise for these graphs.
\end{example}

Table~\ref{tab:six} compares $K_{3,3}$ with three other triangle-free graphs on six vertices.
The values can be computed directly from the defining rules.
The formulas for paths, cycles and stars in~\cite{Choi-Park2015} give the same values.
Although $K_{1,5}$ has fewer edges than $C_6$, its second and third $a$-numbers are larger.
Thus the number of edges alone does not determine the ordering of the higher $a$-numbers.

\begin{table}
\centering
\label{tab:six}
\begin{tabular}{lrrr}
\toprule
$G$ & $a_1(G)$ & $a_2(G)$ & $a_3(G)$ \\
\midrule
$P_6$       & 5 & 9  & 5  \\
$C_6$       & 6 & 15 & 10 \\
$K_{1,5}$   & 5 & 20 & 16 \\
$T_2(6)=K_{3,3}$ & 9 & 39 & 31 \\
\bottomrule
\end{tabular}
\caption{Graph $a$-numbers}
\end{table}

We now compute the $a$-numbers of complete multipartite graphs.
For $K=K_{n_1,\dots,n_r}$, put $n=n_1+\cdots+n_r$ and
$$
B_{2j}=\sum_{\ell=1}^r\binom{n_\ell}{2j}\qquad(j\ge1).
$$
We use the convention $\binom{u}{v}=0$ for integers $v>u\ge0$.

\begin{lemma}\label{lem:multipartite-values}
If $0\le2i\le n$, then
\begin{equation}\label{eq:multipartite-closed}
a_i(K_{n_1,\dots,n_r})
=E_{2i}\binom{n}{2i}
+\sum_{j=1}^i(-1)^jE_{2i-2j}
\binom{n-2j}{2i-2j}B_{2j}.
\end{equation}
In particular, $a_i(K_n)=E_{2i}\binom{n}{2i}$.
\end{lemma}

\begin{proof}
Every nonempty independent set in a complete multipartite graph lies in a unique part.
Hence $I_{2j}(K)=B_{2j}$ for $j\ge1$.
Apply Theorem~\ref{prop:independent} and the equality $a_i(K)=\ha_i(K)$ from Proposition~\ref{prop:compare}.
\end{proof}

The complete multipartite case of this formula can also be obtained from the generating function stated in~\cite{Seo-Shin2015}
$$
\sum_{k_1,\dots,k_r\ge0}a(K_{k_1,\dots,k_r})
\prod_{\ell=1}^r\frac{x_\ell^{k_\ell}}{k_\ell!}
=\frac{1-r+\sum_{\ell=1}^r\cos x_\ell}{\cos(x_1+\cdots+x_r)}.
$$
Thus Lemma~\ref{lem:multipartite-values} is also a finite coefficient form of that known generating function.

For the first three degrees,
\begin{align*}
a_1(K)&=\binom n2-B_2,\\
a_2(K)&=5\binom n4-\binom{n-2}{2}B_2+B_4,\\
a_3(K)&=61\binom n6-5\binom{n-2}{4}B_2
+\binom{n-4}{2}B_4-B_6.
\end{align*}
The identity for $a_i(K)$ is valid for $n\ge2i$.

\begin{lemma}\label{cor:turan-values}
Write $n=qr+s$, where $0\le s<r$, and put
$$
B_{2j}(n,r)=(r-s)\binom q{2j}+s\binom{q+1}{2j}.
$$
For $0\le2i\le n$, we have
\begin{equation}\label{eq:turan-closed}
a_i(T_r(n))=E_{2i}\binom n{2i}
+\sum_{j=1}^i(-1)^jE_{2i-2j}
\binom{n-2j}{2i-2j}B_{2j}(n,r).
\end{equation}
\end{lemma}

\begin{proof}
The graph~$T_r(n)$ has $r-s$ parts of size~$q$ and $s$ parts of
size~$q+1$.
Applying Lemma~\ref{lem:multipartite-values} gives the formula.
\end{proof}

Table~\ref{tab:turan} lists values of \eqref{eq:turan-closed}.

\begin{table}
\centering
\label{tab:turan}
\begin{tabular}{rrrrrrr}
\toprule
$n$ & $r$ & $a_1$ & $a_2$ & $a_3$ & $a_4$ & $a_5$ \\
\midrule
4 & 2 & 4 & 3 & 0 & 0 & 0 \\
4 & 3 & 5 & 4 & 0 & 0 & 0 \\
4 & 4 & 6 & 5 & 0 & 0 & 0 \\
\midrule
5 & 2 & 6 & 13 & 0 & 0 & 0 \\
5 & 3 & 8 & 19 & 0 & 0 & 0 \\
5 & 4 & 9 & 22 & 0 & 0 & 0 \\
5 & 5 & 10 & 25 & 0 & 0 & 0 \\
\midrule
6 & 2 & 9 & 39 & 31 & 0 & 0 \\
6 & 3 & 12 & 57 & 46 & 0 & 0 \\
6 & 4 & 13 & 63 & 51 & 0 & 0 \\
6 & 6 & 15 & 75 & 61 & 0 & 0 \\
\midrule
8 & 2 & 16 & 172 & 820 & 663 & 0 \\
8 & 3 & 21 & 245 & 1183 & 958 & 0 \\
8 & 4 & 24 & 290 & 1408 & 1141 & 0 \\
8 & 8 & 28 & 350 & 1708 & 1385 & 0 \\
\midrule
10 & 2 & 25 & 500 & 5960 & 28915 & 23431 \\
10 & 3 & 33 & 715 & 8625 & 41904 & 33962 \\
10 & 4 & 37 & 826 & 10010 & 48661 & 39441 \\
10 & 10 & 45 & 1050 & 12810 & 62325 & 50521 \\
\bottomrule
\end{tabular}
\caption{Values of $a_i(T_r(n))$.}
\end{table}

Throughout this section, the degree $i\ge1$ is fixed as the number of vertices tends to infinity.
For a graph $G$ on $n\ge2i$ vertices, define
$$
d_i(G)=\frac{a_i(G)}{E_{2i}\binom n{2i}}
=\frac{a_i(G)}{a_i(K_n)}.
$$
By Theorem~\ref{thm:main}, $0\le d_i(G)\le1$.
For~$i=1$, this is the usual edge density.

Let~$H$ be a fixed graph with at least one edge.
If~$n\ge m\ge2i$ and~$S$ is a uniformly chosen $m$-subset
of~$V(G)$, then
\begin{equation}\label{eq:density-average}
\mathbb E\bigl[d_i(G[S])\bigr]=d_i(G).
\end{equation}
Indeed, each $2i$-subset of~$V(G)$ is contained in the same number
of $m$-subsets.
Since every induced subgraph of an $H$-free graph is $H$-free,
the sequence
$$
\frac{\ex_i(n,F)}
{E_{2i}\binom n{2i}}
\qquad (n\ge2i)
$$
is nonincreasing and bounded below by zero.
We therefore define the \emph{$i$th $a$-Tur\'an density of~$H$} by
\begin{equation}\label{eq:density-definition}
\pi_i(H)
=
\lim_{n\to\infty}
\frac{\ex_i(n,F)}
{E_{2i}\binom n{2i}}.
\end{equation}

For~$r,i\ge1$, put
\begin{equation}\label{eq:turan-density}
\tau_i(r)
=
1+\frac1{E_{2i}}
\sum_{j=1}^i
(-1)^j
\binom{2i}{2j}
E_{2i-2j}r^{1-2j}.
\end{equation}
The following proposition shows that~$\tau_i(r)$ is the limiting
density of the Tur\'an graphs.

\begin{lemma}\label{prop:turan-density}
For fixed $r,i\ge1$,
$$
a_i(T_r(n))=E_{2i}\tau_i(r)\binom n{2i}+O_{r,i}(n^{2i-1}).
$$
In particular,
$$
\lim_{n\to\infty}d_i(T_r(n))=\tau_i(r).
$$
\end{lemma}

\begin{proof}
Each part of~$T_r(n)$ has size $n/r+O(1)$.
Hence, for fixed~$j\geq1$,
$$
B_{2j}(n,r)
=
\frac{r^{1-2j}}{(2j)!}n^{2j}
+
O_{r,j}(n^{2j-1}).
$$
It follows that
$$
\frac{\binom{n-2j}{2i-2j}B_{2j}(n,r)}
{\binom{n}{2i}}
=
\binom{2i}{2j}r^{1-2j}
+
O_{r,i}(n^{-1}).
$$
By~\eqref{eq:turan-closed},
\begin{align*}
\frac{a_i(T_r(n))}{\binom{n}{2i}}
&=
E_{2i}
+
\sum_{j=1}^i
(-1)^jE_{2i-2j}
\frac{\binom{n-2j}{2i-2j}B_{2j}(n,r)}
{\binom{n}{2i}}\\
&=
E_{2i}
+
\sum_{j=1}^i
(-1)^jE_{2i-2j}
\binom{2i}{2j}r^{1-2j}
+
O_{r,i}(n^{-1})\\
&=
E_{2i}\tau_i(r)+O_{r,i}(n^{-1}).
\end{align*}
Multiplying by~$\binom{n}{2i}$ proves the lemma.
\end{proof}

We use a standard consequence of the clique blow-up argument.
If $H$ is fixed and $\chi(H)=q\ge2$, then every $H$-free graph on~$n$ vertices has~$o(n^q)$ copies of~$K_q$, uniformly over all such graphs.
Indeed, $H$ is a subgraph of a fixed blow-up of~$K_q$, so this is a special case of~\cite[Proposition~2.1]{Alon-Shikhelman2016}.
Isolated vertices of $H$ do not affect this statement, since they can be added to any copy of the remaining graph once~$n\ge|V(H)|$.

\begin{proof}[Proof of Theorem~\ref{thm:erdos-stone}]
Put $r=\chi(H)-1$.
The graph~$T_r(n)$ is $r$-partite and hence~$H$-free, which gives the lower bound.

For the upper bound, fix an integer $m\ge\max\{2i,r+1\}$.
Let $G$ be any $H$-free graph on $n\ge m$ vertices, and choose an $m$-subset $S$ uniformly at random.
By the preceding clique bound, $G$ has $o(n^{r+1})$ copies of $K_{r+1}$.
Each fixed copy is contained in $S$ with probability
$$
\frac{\binom{n-r-1}{m-r-1}}{\binom nm}=O_m(n^{-r-1}).
$$
Thus the probability $p_n$ that $G[S]$ contains a $K_{r+1}$ is $o(1)$, uniformly over all $H$-free graphs~$G$.
On the complementary event, Theorem~\ref{thm:main} gives
$$
d_i(G[S])\le d_i(T_r(m)).
$$
On the exceptional event, we use $d_i(G[S])\le1$.
By \eqref{eq:density-average},
$$
d_i(G)\le(1-p_n)d_i(T_r(m))+p_n
\le d_i(T_r(m))+o(1).
$$
First let $n\to\infty$, keeping $m$ fixed.
Then let $m\to\infty$ and apply Lemma~\ref{prop:turan-density}.
This gives $\pi_i(H)\le\tau_i(r)$.
The lower bound gives the reverse inequality, and Lemma~\ref{prop:turan-density} yields the theorem.
\end{proof}

For every bipartite $H$ with an edge, $\pi_i(H)=0$ and hence $\ex_i(n,F)=o(n^{2i})$ for each fixed $i\ge1$.
By \eqref{eq:betti}, the theorem also determines the leading asymptotic maximum of $\beta_i(M(G))$ over all $H$-free graphs on $n$ vertices.

The first three densities are
\begin{align*}
\tau_1(r)&=1-\frac1r,\\
\tau_2(r)&=1-\frac6{5r}+\frac1{5r^3},\\
\tau_3(r)&=1-\frac{75}{61r}+\frac{15}{61r^3}-\frac1{61r^5}.
\end{align*}
Setting $\tau_0(r)=1$, we can also write \eqref{eq:turan-density} as
$$
\sum_{i\ge0}E_{2i}\tau_i(r)\frac{x^{2i}}{(2i)!}
=\frac{1-r+r\cos(x/r)}{\cos x}.
$$
For example, $\tau_1(2)=1/2$, $\tau_2(2)=17/40$ and $\tau_3(2)=811/1952$.
Table~\ref{tab:density} gives further values, rounded to six decimal places.

\begin{table}[htbp]
\centering
\label{tab:density}
\begin{tabular}{rrrrr}
\toprule
$r$ & $\tau_1(r)$ & $\tau_2(r)$ & $\tau_3(r)$ & $\tau_4(r)$ \\
\midrule
2 & 0.500000 & 0.425000 & 0.415471 & 0.414356 \\
3 & 0.666667 & 0.607407 & 0.599204 & 0.598206 \\
4 & 0.750000 & 0.703125 & 0.696449 & 0.695626 \\
5 & 0.800000 & 0.761600 & 0.756060 & 0.755373 \\
10 & 0.900000 & 0.880200 & 0.877295 & 0.876931 \\
\bottomrule
\end{tabular}
\caption{Normalized $a$-Tur\'an densities $\tau_i(r)$.}
\end{table}

\begin{remark}
For fixed~$r\ge2$, the simple poles at $x=\pm\pi/2$ and $x=\pm3\pi/2$ in the preceding generating function give~$\tau_i(r)=1-r+r\cos(\pi/(2r))+O_r(9^{-i})$ as~$i\to\infty$, which explains the rapid stabilization in Table~\ref{tab:density}.
\end{remark}

\end{document}